\documentclass[12pt]{article}
\usepackage{amsmath,amssymb}
\usepackage{theorem}

\newtheorem{theorem}{Theorem}[section]
\newtheorem{proposition}[theorem]{Proposition}
\newtheorem{lemma}[theorem]{Lemma}
\newtheorem{pfpf}{{\it Proof.}}

\newenvironment{prf}{\begin{pfpf}\rm}{\hspace*{\fill}{$\square$}\end{pfpf}}
\theorembodyfont{\rmfamily}\newtheorem{remark}[theorem]{Remark}

\numberwithin{equation}{section}

\newcommand{\Td}{\operatorname{Td}}
\newcommand{\rootht}{\operatorname{ht}}

\begin{document}

\title{\textbf{Power sums of {C}oxeter exponents}\\[9mm]}

\author{\Large John M.\ Burns\\[3mm]
{\textit{School of Mathematics, Statistics and Applied Mathematics}}\\
{\textit{National University of Ireland, Galway}}\\
{\textit{University Road, Galway, Ireland}}\\[2mm]
e-mail: \texttt{john.burns@nuigalway.ie}\\[5mm]
\and
\Large Ruedi Suter\\[3mm]
{\textit{Departement Mathematik}}\\
{\textit{ETH Z\"urich}}\\
{\textit{R\"amistrasse 101, 8092 Z\"urich, Switzerland}}\\[2mm]
e-mail: \texttt{suter@math.ethz.ch}\\[5mm]
}

\maketitle

\begin{abstract}
Consider an irreducible finite Coxeter system. We show that for any
nonnegative integer $n$ the sum of the $n$th powers of the Coxeter
exponents can be written uniformly as a polynomial in four parameters:
$h$ (the Coxeter number), $r$ (the rank), $\alpha$, $\beta$ (two
further parameters).
\end{abstract}

\section{Introduction}
Let $(W,S)$ be an irreducible finite Coxeter system of rank $r$ with
$S=\{s_1,\dots,s_r\}$ its set of simple reflections. The Coxeter transformation
$c:=s_1\ldots s_r\in W$ has order $|c|=h$ known as the Coxeter number,
and the eigenvalues of $c$ in the reflection representation of $W$ are
of the form $e^{2\pi i m_1/h},\dots,e^{2\pi i m_r/h}$ with
$1=m_1\leqslant m_2\leqslant\cdots\leqslant m_r=h-1$ the exponents of
$(W,S)$. Furthermore, for any permutation $\sigma$ of $\{1,\dots,r\}$
the elements $c$ and $s_{\sigma(1)}\ldots s_{\sigma(r)}$ are conjugate
in $W$. Hence the exponents do not depend on the enumeration of the
simple reflections. Recall that the symmetry $m_i+m_{r+1-i}=h$ follows
from the facts that $c$ has no eigenvalue $1$ and that the reflection
representation is defined over the reals.

In this note we will derive uniform expressions for the power sums
$\sum_{i=1}^r m_i^{\,n}$ for any $n\in\mathbb Z_{\geqslant0}$. Of course, for $n=0$
the sum is $r$, and for $n=1$ the symmetry $m_i+m_{r+1-i}=h$ shows that
the sum is $\frac12 rh$. We shall see that
$$\sum_{i=1}^r m_i^{\,n}=n!\,r\Td_n(\gamma_1,\dots,\gamma_n)$$
where $\Td_n(\gamma_1,\dots,\gamma_n)$ denotes the $n$th Todd polynomial
evaluated at $\gamma_1,\dots,\gamma_n$ (for $n$ odd
$\Td_n(\gamma_1,\dots,\gamma_n)$ does not depend on $\gamma_n$, as follows
from Proposition~\ref{ToddSymm}).
The $\gamma_i$'s can be chosen
to be polynomials in four parameters (details below) with integer
coefficients. 
This answers Panyushev's question in \cite{Pa}.

\section{Some history and preliminaries}
For type $\mathsf A_r$ the exponents are just $1,2,\dots,r$ and one has
Bernoulli's formula
\begin{equation}\label{Bernoulli}
\sum_{i=1}^r i^n=\frac{1}{n+1}\bigl(B_{n+1}(r+1)-B_{n+1}(1)\bigr)
\end{equation}
where $B_{n+1}(x)$ is the $(n+1)$st Bernoulli polynomial,
defined by the expansion
$$\sum_{n=0}^\infty B_n(x)\frac{t^n}{n!}=\frac{t\,e^{tx}}{e^t-1}\,.$$

For general types uniform formulae for the power sums up
to third power are listed in the epilogue of \cite{Su}.
Besides the Coxeter number $h$ and the rank $r$
they depend (for the squares and the cubes) on a further parameter
$\gamma$ which is defined for the crystallographic types with
crystallographic root system $\Phi$ (${}=\Phi_+\cup\Phi_-$ a decomposition
into the sets of positive and negative roots) by the formula
(see \cite[Ch.~VI, \S~1, no.~12]{Bou})
\begin{equation}\label{definitiongamma}
\sum_{\varphi\in\Phi}\frac{\langle\lambda|\varphi\rangle
\langle\mu|\varphi\rangle}{\langle\varphi|\varphi\rangle^2}=
\gamma\,\langle\lambda|\mu\rangle\qquad(\lambda,\mu\in
\operatorname{span}_\mathbb R\Phi)
\end{equation}
where $\langle\phantom{\varphi}|\phantom{\varphi}\rangle$
denotes the Killing form on $\operatorname{span}_\mathbb R\Phi$, which
is the $W$-invariant (symmetric) bilinear form characterized by
$$\langle\lambda|\mu\rangle=\sum_{\varphi\in\Phi}
\langle\lambda|\varphi\rangle\langle\mu|\varphi\rangle
\qquad(\lambda,\mu\in\operatorname{span}_\mathbb R\Phi).$$
It turns out that $\gamma=kgg^\vee$ where $k=\langle\theta|\theta\rangle
/\langle\theta_{\textrm{s}}|\theta_{\textrm{s}}\rangle\in\{1,2,3\}$
with $\theta,\theta_{\textrm{s}}\in\Phi_+$ the highest resp.\ highest
short roots, and $g=1/\langle\theta|\theta\rangle
\in\mathbb Z_{>0}$ is the dual Coxeter number of $\Phi$ whereas $g^\vee$ is the
dual Coxeter number of the dual root system $\Phi^\vee$. So $\gamma=h^2$
if $\Phi$ is simply-laced. 
For the noncrystallographic types $\gamma=2m^2-5m+6$ for $\mathsf I_2(m)$
(the formula is also valid for the crystallographic types, where
$m=3,4,6$); $\gamma=124$ for type $\mathsf H_3$; and $\gamma=1116$ for
type $\mathsf H_4$.

The formulae from \cite{Su} read as follows:
\begin{equation}\label{lowpowers}
\sum_{i=1}^r m_i^{\,n}=\begin{cases}
r&\mbox{if $n=0$,}\\[2mm]
\tfrac12rh&\mbox{if $n=1$,}\\[2mm]
\tfrac16r(h^2+\gamma-h)&\mbox{if $n=2$,}\\[2mm]
\tfrac14rh(\gamma-h)&\mbox{if $n=3$.}
\end{cases}
\end{equation}

\begin{remark}
The power sum for the fourth powers
is not of the form $r$ times a function depending
only on $h$ and $\gamma$, as a computation
for the types $\mathsf A_{h-1}$ and $\mathsf D_{(h+2)/2}$ shows.
\end{remark}

Panyushev recently gave the universal formula
\cite[Proposition~3.1]{Pa}
\begin{equation}\label{heightssquaresum}
\sum_{\varphi\in\Phi_+}\rootht(\varphi)^2=\frac{1}{12}r(h+1)\gamma
\end{equation}
for the sum of the heights squares of all positive roots. He then
suspects \cite[Remark~3.4]{Pa} that for the sum of the heights of
all positive roots there is no similar formula in the general case;
however, for simply-laced root systems he mentions
\begin{equation}\label{heightssum}
\sum_{\varphi\in\Phi_+}\rootht(\varphi)=\frac16r(h^2+h)
\end{equation}
and asks for which values of $n$ there is a nice closed expression for
$\sum_{\varphi\in\Phi_+}\rootht(\varphi)^n$.
Our result shows that there are universal formulae for all
$n\in\mathbb Z_{\geqslant0}$.
In fact, let $(k_1,\dots,k_{h-1})$ be the partition dual to $(m_r,\dots,m_1)$;
then it is well-known (see, e.\,g., \cite[Section~3.20]{Hu})
that there are exactly $k_j$ roots of height $j$
in $\Phi_+$. Hence
\begin{equation}\label{heightspowersum}
\sum_{\varphi\in\Phi_+}\rootht(\varphi)^n=\sum_{i=1}^r\bigl(1^n+2^n+\cdots
+m_i^{\,n}\bigr).
\end{equation}
In particular, using (\ref{lowpowers}) we recover (\ref{heightssquaresum})
and have
\begin{equation}\label{generalsum}
\sum_{\varphi\in\Phi_+}\rootht(\varphi)=
\sum_{i=1}^r\frac{m_i^{\,2}+m_i}{2}=\frac{1}{12}r\bigl(h^2+\gamma+2h\bigr)
\end{equation}
which generalizes (\ref{heightssum}) to all types.

Alternatively, using the symmetry $m_i+m_{r+1-i}=h$ we can write as in
\cite[Proposition 2.1]{Bu}
\begin{equation}\label{symmetrycubes}
h^2\sum_{i=1}^r m_i-3h\sum_{i=1}^r m_i^{\,2}+2\sum_{i=1}^rm_i^{\,3}=0.
\end{equation}
Hence
\begin{align*}
\sum_{\varphi\in\Phi_+}\rootht(\varphi)^2&\stackrel{(\ref{heightspowersum})}{=}
\sum_{i=1}^r\frac{m_i(m_i+1)(2m_i+1)}{6}
=\sum_{i=1}^r\frac{m_i^{\,3}}{3}+\sum_{i=1}^r\frac{m_i^{\,2}}{2}
+\sum_{i=1}^r\frac{m_i}{6}\\
&\stackrel{(\ref{symmetrycubes})}{=}
-h^2\sum_{i=1}^r\frac{m_i}{6}+h\sum_{i=1}^r\frac{m_i^{\,2}}{2}
+\sum_{i=1}^r\frac{m_i^{\,2}}{2}
+\sum_{i=1}^r\frac{m_i}{6}\\
&\stackrel{\phantom{(0.0)}}{=}(h+1)\sum_{i=1}^r\frac{m_i(m_i+1)}{2}
-\Bigl(\frac{h+1}{2}+\frac{h^2-1}{6}\Bigr)
\raisebox{0pt}[0pt][0pt]{$\displaystyle
\underbrace{\sum_{i=1}^r\ m_i}_{\displaystyle{}=\frac{rh}{2}}$}\\
&\stackrel{(\ref{heightspowersum})}{=}
(h+1)\sum_{\varphi\in\Phi_+}\rootht(\varphi)-(h+1)\frac{rh(h+2)}{12}
\end{align*}
so that (\ref{generalsum}) is recovered from (\ref{heightssquaresum}).

We shall stick to the exponents rather than the heights in order not to
restrict our considerations to the crystallographic types.

\section{Power sums and Todd polynomials}
\newlength{\lenga}
\settowidth{\lenga}{$\tfrac16r(h^2+\gamma-h)$}
The observation that (\ref{lowpowers}) can be written as
\begin{equation}\label{smallpowers}
\sum_{i=1}^r m_i^{\,n}=\begin{cases}
\makebox[\lenga][l]{$r$}=0!\,r\Td_0&\mbox{if $n=0$,}\\[2mm]
\makebox[\lenga][l]{$\tfrac12rh$}=1!\,r\Td_1(h)&\mbox{if $n=1$,}\\[2mm]
\makebox[\lenga][l]{$\tfrac16r(h^2+\gamma-h)$}=2!\,r\Td_2(h,\gamma-h)
&\mbox{if $n=2$,}\\[2mm]
\makebox[\lenga][l]{$\tfrac14rh(\gamma-h)$}=3!\,r\Td_3(h,\gamma-h,\ast)
&\mbox{if $n=3$,}
\end{cases}
\end{equation}
where $\Td_0=1$, $\Td_1(c_1)=\frac12 c_1$, $\Td_2(c_1,c_2)=\frac{1}{12}
(c_1^{\,2}+c_2)$, and $\Td_3(c_1,c_2,c_3)=\frac{1}{24}c_1c_2$ are Todd polynomials
(the general definition will be recalled in the proof of Theorem~\ref{mainthm}),
suggests the ansatz
\begin{equation}\label{ansatz}
\sum_{i=1}^r m_i^{\,n}=n!\,r\Td_n(\gamma_1,\dots,\gamma_n).
\end{equation}
From (\ref{smallpowers}) and (\ref{ansatz}) we get
\begin{equation}\label{gammaonetwo}
\mbox{$\gamma_1=h$ and $\gamma_2=\gamma-h$}
\end{equation}
and are looking for
solutions $\gamma_3,\gamma_4,\ldots$. Note that the symmetry
$m_i+m_{r+1-i}=h$ implies the identities (for $a,b\in\mathbb Z_{\geqslant0}$)
\begin{equation}\label{symmetry}
\sum_{j=0}^a(-1)^{a-j}\binom{a}{j}h^j\sum_{i=1}^r m_i^{\,a+b-j}=
\sum_{j=0}^b(-1)^{b-j}\binom{b}{j}h^j\sum_{i=1}^r m_i^{\,a+b-j}
\end{equation}
that generalize (\ref{symmetrycubes}), which is (\ref{symmetry}) for
$\{a,b\}=\{1,2\}$.

\begin{proposition}\label{ToddSymm}
For $a,b\in\mathbb Z_{\geqslant0}$ one has the identity
\begin{align}\notag
&\sum_{j=0}^a(-1)^{a-j}\binom{a}{j}\,c_1^j(a+b-j)!
\Td_{a+b-j}(c_1,\dots,c_{a+b-j})\\\label{lemma}
&\qquad=
\sum_{j=0}^b(-1)^{b-j}\binom{b}{j}\,c_1^j(a+b-j)!
\Td_{a+b-j}(c_1,\dots,c_{a+b-j}).
\end{align}
\end{proposition}
\begin{prf}
For instance, one verifies the formula (\ref{lemma}) for $a=0$ and
all $b\in\mathbb Z_{\geqslant0}$ by using a generating series and then proceeds
by induction on $a$.
\end{prf}

Strictly speaking we don't need Proposition~\ref{ToddSymm}. But it is worth
noting that it indicates that we seem to be on the right track when using
the ansatz~(\ref{ansatz}).

\begin{lemma}
Let $m_1\leqslant\dots\leqslant m_r\in\mathbb Z_{>0}$ be such that there are
multisets $V_+$ and $V_-$ of positive integers satisfying
\begin{equation}\label{expsum}
\sum_{i=1}^r q^{m_i}=\frac{q\prod_{v\in V_+}(1-q^v)}{\prod_{v\in V_-}(1-q^v)}\,.
\end{equation}
Then
\begin{align}\label{ABrab}
\prod_{v\in V_+}v&=r\prod_{v\in V_-}v\\\label{cardinality}
|V_+|&=|V_-|.
\end{align}
\end{lemma}
\begin{prf}
The equality (\ref{ABrab}) is clear from the $q\to1$ limit in (\ref{expsum});
(\ref{cardinality}) follows since $1-q^v$ has exactly one factor $1-q$ and
the polynomial on the left hand side in (\ref{expsum}) has
neither a zero nor a pole at $q=1$. Note also that $m_1=1$ and
$m_2>1$ if $r\geqslant2$.
\end{prf}

\begin{theorem}\label{mainthm}
Let $m_1\leqslant\dots\leqslant m_r\in\mathbb Z_{>0}$ be such that there are
multisets $V_+$ and $V_-$ of positive integers satisfying
\begin{equation*}\tag{\ref{expsum}}
\sum_{i=1}^r q^{m_i}=\frac{q\prod_{v\in V_+}(1-q^v)}{\prod_{v\in V_-}(1-q^v)}\,.
\end{equation*}
We fix (for simplicity) a positive integer $p$ and define
$\gamma_0({}=1),\gamma_1,\gamma_2,\gamma_3,\ldots$ by the generating series
\begin{equation}\label{defgamma}
\sum_{n=0}^\infty\gamma_n t^n=\frac{\prod_{v\in V_-}(1-vt)}{\prod_{v\in V_+}
(1-vt)}\sqrt[\raisebox{0.3em}{$\textstyle p$}]{\frac{1+pt}{1-pt}}\,.
\end{equation}
Then for $n\in\mathbb Z_{\geqslant0}$
\begin{equation}\label{formula}
\sum_{i=1}^r m_i^{\,n}=n!\,r\Td_n(\gamma_1,\dots,\gamma_n).
\end{equation}
\end{theorem}
\begin{prf}
We consider the exponential generating series (with $q:=e^t$) of both
sides in (\ref{formula})
\begin{align}\label{LHS}
\sum_{n=0}^\infty\Bigl(\sum_{i=1}^r m_i^{\,n}\Bigr)\frac{t^n}{n!}
&=\sum_{i=1}^r e^{m_i t}=\sum_{i=1}^r q^{m_i}
\stackrel{(\ref{expsum})}{=}
\frac{q\prod_{v\in V_+}(1-q^v)}{\prod_{v\in V_-}(1-q^v)}\\\label{RHS}
\sum_{n=0}^\infty\Bigl(n!\,r\Td_n(\gamma_1,\dots,\gamma_n)\Bigr)\frac{t^n}{n!}
&=r\sum_{n=0}^\infty \Td_n(\gamma_1,\dots,\gamma_n)t^n
=r\prod_{j=1}^\infty\frac{x_jt}{1-e^{-x_jt}}
\end{align}
where the last equality incorporates the definition of the Todd polynomials
by means of their generating series in $t$ with coefficients in the elementary
symmetric functions in $x_1,x_2,\ldots$ so that 
\begin{equation*}
\prod_{j=1}^\infty(1+x_jt)=\sum_{n=0}^\infty\gamma_n t^n
\end{equation*}
and hence by (\ref{defgamma})
\begin{equation*}
\frac{(1-pt)\prod_{v\in V_+}(1-vt)^p}{(1+pt)\prod_{v\in V_-}(1-vt)^p}
\prod_{j=1}^\infty(1+x_jt)^p=1
\end{equation*}
(that is, the supersymmetric elementary symmetric functions
in $-p$, $-v$ ($p$ times, for every $v\in V_+$), $x_1$ ($p$ times),
$x_2$ ($p$ times), \ldots;
$-p$, $v$ ($p$ times, for every $v\in V_-$) all vanish) so that the formal
expansion
\begin{equation*}
\underbrace{\Bigl(\frac{-pt}{1-e^{pt}}\Bigr)\Bigl(
\frac{1-e^{-pt}}{pt}\Bigr)}_{\textstyle{}=e^{-pt}}
\prod_{v\in V_+}\Bigl(\frac{-vt}{1-e^{vt}}\Bigr)^p
\prod_{v\in V_-}\Bigl(\frac{1-e^{vt}}{-vt}\Bigr)^p
\Bigl(\prod_{j=1}^\infty\frac{x_jt}{1-e^{-x_jt}}\Bigr)^p=1
\end{equation*}
or taking $p$th roots (look at $t=0$ to choose the correct branch) 
\begin{equation*}\prod_{j=1}^\infty\frac{x_jt}{1-e^{-x_jt}}=e^t
\prod_{v\in V_+}\Bigl(\frac{1-e^{vt}}{-vt}\Bigr)
\prod_{v\in V_-}\Bigl(\frac{-vt}{1-e^{vt}}\Bigr).
\end{equation*}
Therefore
we can write the right hand side in (\ref{RHS}) as (recall $q=e^t$)
\begin{equation*}
r\prod_{j=1}^\infty\frac{x_jt}{1-e^{-x_jt}}=
\underbrace{\frac{r\prod_{v\in V_-}v}{\prod_{v\in V_+}v}}_{\textstyle{}=1}\cdot
\frac{q\prod_{v\in V_+}(1-q^v)}{\prod_{v\in V_-}(1-q^v)}=
\frac{q\prod_{v\in V_+}(1-q^v)}{\prod_{v\in V_-}(1-q^v)}
\end{equation*}
where we have used (\ref{cardinality}) $|V_+|=|V_-|$ to cancel factors $t$
and then (\ref{ABrab}) to simplify the product.
Thus the right hand side of (\ref{RHS})
is identical to the right hand side of (\ref{LHS}), which proves
(\ref{formula}).
\end{prf}

\begin{remark}
Instead of the definition (\ref{defgamma}) for
$\gamma_0,\gamma_1,\gamma_2,\gamma_3,\ldots$ one can define more generally
\begin{equation*}
\sum_{n=0}^\infty\gamma_n t^n=\frac{\prod_{v\in V_-}(1-vt)}{\prod_{v\in V_+}(1-vt)}
\prod_{k=1}^K\Bigl(\frac{1+\pi_kt}{1-\pi_kt}\Bigr)^{\mu_k}
\end{equation*}
with $\pi_1,\dots,\pi_K\in\mathbb R$ and $\mu_1,\dots,\mu_K\in\mathbb Q$
satisfying $\sum_{k=1}^K\pi_k\mu_k=1$ (and 
for general $m_1$ (with $q^{m_1}$ instead of
$q$ as first factor in the right hand side of (\ref{expsum}))
just require that
$\sum_{k=1}^K\pi_k\mu_k=m_1$).
\end{remark}

\section{Root system considerations}
To apply Theorem~\ref{mainthm} in the context of root systems we need the
following proposition.

\begin{proposition}\label{rootsystcyclotomic}
Let $m_1\leqslant\cdots\leqslant m_r$ be the exponents of an irreducible
(crystallographic (and reduced) or noncrystallographic) finite
root system (of rank $r$). Then there are multisets $V_+$ and $V_-$ of
positive integers such that
\begin{equation}\tag{\ref{expsum}}
\sum_{i=1}^r q^{m_i}=\frac{q\prod_{v\in V_+}(1-q^v)}{\prod_{v\in V_-}(1-q^v)}\,.
\end{equation}
Furthermore,
$|V_\pm|\leqslant2$ if \,$V_+\cap V_-=\varnothing$.
\end{proposition}
\begin{prf}
According to the first note added in proof in \cite{Sa} I.~G.~Macdonald
was acquainted with the fact that (\ref{expsum}) holds for all irreducible
finite Coxeter groups.

The classification shows that the following three cases exhaust all
possible types.
\begin{itemize}
\item[(1)] For the types $\mathsf A_r$, $\mathsf C_r/\mathsf B_r$, and
types of rank ${}\leqslant3$ the sequence of exponents forms an
arithmetic progression $1,m_2,\dots,1+(r-1)(m_2-1)$ (or just $1$ if $r=1$).
Hence
$$\sum_{i=1}^r q^{m_i}=
\begin{cases}q&\mbox{if $r=1$}\\
\dfrac{q(1-q^{r(m_2-1)})}{1-q^{m_2-1}}&\mbox{if $r\geqslant 2$}
\end{cases}$$
so that we can take $V_+=V_-=\varnothing$ if $r=1$ and
$V_+=\{r(m_2-1)\}$ and $V_-=\{m_2-1\}$ if $r\geqslant2$.
\item[(2)] For the types of rank $4$ we have
$$\sum_{i=1}^4 q^{m_i}=q+q^{m_2}+q^{h-m_2}+q^{h-1}
=\frac{q(1-q^{2(m_2-1)})(1-q^{2(h-m_2-1)})}{(1-q^{m_2-1})(1-q^{h-m_2-1})}$$
so that we can take $V_+=\{2(m_2-1),2(h-m_2-1)\}$ and $V_-=\{m_2-1,h-m_2-1\}$.
\item[(3)] For the simply-laced types ($\mathsf{ADE}$) the root system
is the Weyl group orbit of the highest root: $\Phi=W\theta$. The stabilizer
of $\theta$ is $W_{\perp\theta}$, the reflection group generated by those
simple reflections in $W$ that fix $\theta$. The root system is thus
isomorphic as a $W$-set to $W/W_{\perp\theta}$. We need the usual length
function $\ell:W\to\mathbb Z_{\geqslant0}$ defined as $\ell(w)=k$ if $w$ can
be written as a product of $k$ but not less than $k$ simple reflections.
If $\varphi=w\theta$
is any positive root with $w$ chosen such that $\ell(w)$ is minimal, then
$\rootht(\varphi)=\rootht(\theta)-\ell(w) =h-1-\ell(w)$. Since the reflection
along a simple root $\psi$ maps
$\psi$ (of height $1$) to $-\psi$ (of height $-1$), we have similarly
the equality
$\rootht(\varphi)=\rootht(\theta)-\ell(w)-1=h-2-\ell(w)$ if $\varphi=w\theta$
is any negative root with $w$ chosen such that $\ell(w)$ is minimal. So we have
\begin{align*}
\sum_{\substack{wW_{\perp\theta}\in W/W_{\perp\theta}\\
\ell(w)\textup{\ minimal}}}q^{\ell(w)}&=
\sum_{\varphi\in\Phi_+}\bigl(q^{h-1-\rootht(\varphi)}+q^{h-2+\rootht(\varphi)}
\bigr)
\intertext{and since $1,\dots,m_1,\,1,\dots,m_2,\,\dots,\,1,\dots,m_r$
enumerates $\rootht(\varphi)$ as $\varphi$ runs over $\Phi_+$, we can continue}
&=\sum_{i=1}^r\sum_{j=1}^{m_i}\bigl(q^{h-1-j}+q^{h-2+j}\bigr)
\intertext{and using the symmetry $m_i+m_{r+1-i}=h$ we obtain}
&=\sum_{i=1}^r\sum_{j=0}^{h-1}q^{m_i-1+j}
=\Bigl(\sum_{i=1}^r q^{m_i-1}\Bigr)\frac{1-q^h}{1-q}\,.
\intertext{On the other hand by the Chevalley-Solomon identity for
the Poincar\'e series of finite Coxeter groups
(see, e.\,g., \cite[Section~3.15]{Hu}) we have}
\sum_{\substack{wW_{\perp\theta}\in W/W_{\perp\theta}\\
\ell(w)\textup{\ minimal}}}q^{\ell(w)}&=
\Bigl(\prod_{i=1}^r\frac{1-q^{m_i+1}}{1-q}\Bigr)
\Bigl(\prod_{i=1}^s\frac{1-q}{1-q^{\widetilde m_i+1}}\Bigr)
\intertext{where $\widetilde m_1,\dots,\widetilde m_s$ lists the exponents
of all the irreducible components of $W_{\perp\theta}$. Since $m_r+1=h$ we
finally get}
\sum_{i=1}^r q^{m_i}&=\frac{q}{(1-q)^{r-s-1}}\frac{\prod_{i=1}^{r-1}
(1-q^{m_i+1})}{\prod_{i=1}^s(1-q^{\widetilde m_i+1})}
\intertext{and the following table finishes the proof. (We have left out
the types $\mathsf A_r$ which were already dealt with in case~(1).)}
\end{align*}
\end{itemize}
\mbox{}\\[-4mm]
\renewcommand{\arraystretch}{1.2}%
\footnotesize%
$$\begin{array}{|cll|cl|c|c|}\hline
\multicolumn{2}{|c}{\mbox{type $W$}}&\mbox{exponents}+1
&\mbox{type $W_{\perp\theta}$}&\mbox{exponents}+1&
V_+&V_-\\\hline\hline
\mathsf D_r&(r\geqslant4)&2,4,\dots,2r-2,r&
\mathsf A_1+\mathsf D_{r-2}&2,2,4,\dots,2r-6,r-2&
\{r,2r-4\}&\{2,r-2\}\\
\mathsf E_6&&2,5,6,8,9,12&\mathsf A_5&2,3,4,5,6&\{8,9\}&\{3,4\}\\
\mathsf E_7&&2,6,8,10,12,14,18&\mathsf D_6&2,4,6,8,10,6&\{12,14\}&\{4,6\}\\
\mathsf E_8&&2,8,12,14,18,20,24,30&\mathsf E_7&2,6,8,10,12,14,18&
\{20,24\}&\{6,10\}\\\hline
\end{array}$$
\normalsize
Multisets are needed for type $\mathsf D_4$.
\end{prf}

Note that for $r\geqslant2$ (\ref{expsum}) implies that $m_2-1\in V_-$.
Furthermore, for all the crystallographic types except $\mathsf A_1$ and
$\mathsf G_2$, $m_2-1=d$ is the largest coefficient of the highest root
(when written as a linear combination of the simple roots). This observation
extends to the noncrystallographic types $\mathsf H_3$ and $\mathsf H_4$
if we define $d=4$ and $d=10$, respectively, as suggested by the
following folding procedure, $\mathsf D_6\rightsquigarrow\mathsf H_3$
and $\mathsf E_8\rightsquigarrow\mathsf H_4$.
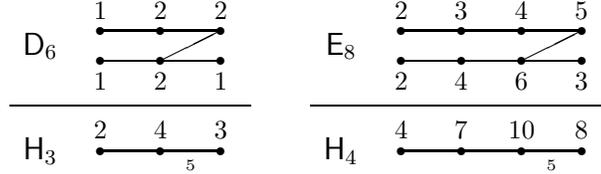
\begin{center}
\setlength{\unitlength}{4mm}
\begin{picture}(20,5)(-4,0)
\put(0,0){\line(1,0){4}}
\put(0,3){\line(1,0){4}}
\put(0,4){\line(1,0){4}}
\put(2,3){\line(2,1){2}}
\put(10,0){\line(1,0){6}}
\put(10,3){\line(1,0){6}}
\put(10,4){\line(1,0){6}}
\put(14,3){\line(2,1){2}}
\multiput(0,0)(2,0){3}{\circle*{0.3}}
\multiput(0,3)(2,0){3}{\circle*{0.3}}
\multiput(0,4)(2,0){3}{\circle*{0.3}}
\multiput(10,0)(2,0){4}{\circle*{0.3}}
\multiput(10,3)(2,0){4}{\circle*{0.3}}
\multiput(10,4)(2,0){4}{\circle*{0.3}}
\put(0,0.4){\makebox(0,0)[b]{\footnotesize$2$}}
\put(2,0.4){\makebox(0,0)[b]{\footnotesize$4$}}
\put(4,0.4){\makebox(0,0)[b]{\footnotesize$3$}}
\put(0,4.4){\makebox(0,0)[b]{\footnotesize$1$}}
\put(2,4.4){\makebox(0,0)[b]{\footnotesize$2$}}
\put(4,4.4){\makebox(0,0)[b]{\footnotesize$2$}}
\put(0,2.6){\makebox(0,0)[t]{\footnotesize$1$}}
\put(2,2.6){\makebox(0,0)[t]{\footnotesize$2$}}
\put(4,2.6){\makebox(0,0)[t]{\footnotesize$1$}}
\put(3,-0.3){\makebox(0,0)[t]{\tiny$5$}}
\put(10,0.4){\makebox(0,0)[b]{\footnotesize$4$}}
\put(12,0.4){\makebox(0,0)[b]{\footnotesize$7$}}
\put(14,0.4){\makebox(0,0)[b]{\footnotesize$10$}}
\put(16,0.4){\makebox(0,0)[b]{\footnotesize$8$}}
\put(10,4.4){\makebox(0,0)[b]{\footnotesize$2$}}
\put(12,4.4){\makebox(0,0)[b]{\footnotesize$3$}}
\put(14,4.4){\makebox(0,0)[b]{\footnotesize$4$}}
\put(16,4.4){\makebox(0,0)[b]{\footnotesize$5$}}
\put(10,2.6){\makebox(0,0)[t]{\footnotesize$2$}}
\put(12,2.6){\makebox(0,0)[t]{\footnotesize$4$}}
\put(14,2.6){\makebox(0,0)[t]{\footnotesize$6$}}
\put(16,2.6){\makebox(0,0)[t]{\footnotesize$3$}}
\put(15,-0.3){\makebox(0,0)[t]{\tiny$5$}}
\put(-3,1.5){\line(1,0){8}}
\put(7,1.5){\line(1,0){10}}
\put(-2,3.5){\makebox(0,0){$\mathsf D_6$}}
\put(-2,0){\makebox(0,0){$\mathsf H_3$}}
\put(8,3.5){\makebox(0,0){$\mathsf E_8$}}
\put(8,0){\makebox(0,0){$\mathsf H_4$}}
\end{picture}
\end{center}
For $\mathsf I_2(m)$ we have $m_2-1=m-2$, but the folding procedure
gives $d=\bigl\lfloor\frac m2\bigr\rfloor$. In fact, for $m=2k+1$ odd
$\mathsf A_{2k}\rightsquigarrow\mathsf I_2(2k+1)$ with $d=k$ (and the other
coefficient is $k$, too). For $m=2k$ even
$\mathsf D_{k+1}\rightsquigarrow\mathsf I_2(2k)$ with $d=k$ (and the other
coefficient is $k-1$); alternatively we can fold $\mathsf A_{2k-1}
\rightsquigarrow\mathsf I_2(2k)$ and also $\mathsf E_6\rightsquigarrow
\mathsf I_2(12)$, $\mathsf E_7\rightsquigarrow\mathsf I_2(18)$, and
$\mathsf E_8\rightsquigarrow\mathsf I_2(30)$.

For $\mathsf A_r$, $\mathsf C_r$, $\mathsf B_r$, $\mathsf I_2(m)$, and
$\mathsf H_3$ one can append the same element(s) to both $V_+$ and $V_-$
to make all the above multisets $V_+$ and $V_-$ have cardinality $2$.

The following proposition gives a uniform description of multisets
$V_+=\{A,B\}$ and $V_-=\{\alpha,\beta\}$ satisfying (\ref{expsum}) 
in terms of three parameters: the Coxeter number $h$, the coefficient $d$, and
$\nu:={}$ the number of times $d$ occurs among the marks in the extended
Dynkin diagram minus $1$, and extended to the noncrystallographic types
as displayed in the following table. The table also shows the values
of $\gamma$ (see (\ref{definitiongamma}) and the text afterwards). Some
parameters $\beta$ (and for type $\mathsf A_1$ also $\alpha$) are irrelevant and
are left unspecified. Clearly, one can interchange $A\leftrightarrow B$
and also $\alpha\leftrightarrow\beta$.
\renewcommand{\arraystretch}{1.1}%
$$\begin{array}{|cl|c|c|c|c|c|c|c|}\hline
\multicolumn{2}{|c|}{\mbox{type}}&r&h&\gamma&d&A,B&\alpha,\beta
&\nu\\\hline\hline
\mathsf A_1&&1&2&4&1&\alpha,\beta&\alpha,\beta&1\\
\mathsf A_r&(r\geqslant2)&r&r+1&(r+1)^2&1&r,\beta&1,\beta&r\\
\mathsf C_r/\mathsf B_r&(r\geqslant2)&r&2r&4r^2+2r-2&2&2r,\beta&2,\beta&r-2\\
\mathsf D_r&(r\geqslant4)&r&2r-2&(2r-2)^2&2&r,2(r-2)&2,r-2&r-4\\
\mathsf E_6&&6&12&144&3&8,9&3,4&0\\
\mathsf E_7&&7&18&324&4&12,14&4,6&0\\
\mathsf E_8&&8&30&900&6&20,24&6,10&0\\
\mathsf F_4&&4&12&162&4&8,12&4,6&0\\
\mathsf G_2\makebox[0pt][l]{${}=\mathsf I_2(6)$}&&2&6&48&3&8,\beta&4,\beta&0\\
\mathsf H_2\makebox[0pt][l]{${}=\mathsf I_2(5)$}&&2&5&31&2&6,\beta&3,\beta&1\\
\mathsf H_3&&3&10&124&4&12,\beta&4,\beta&0\\
\mathsf H_4&&4&30&1116&10&20,36&10,18&0\\
\mathsf I_2(2k+1)&(k\geqslant3)&2&2k+1&8k^2-2k+3&k&4k-2,\beta&2k-1,\beta&1\\
\mathsf I_2(2k)&(k\geqslant4)&2&2k&8k^2-10k+6&k&4k-4,\beta&2k-2,\beta&0\\\hline
\end{array}$$

$$\begin{array}{|cl|c|c|c|c|c|c|c|}\hline
\multicolumn{9}{|l|}{\mbox{redefined parameters $d$ and $\nu$
for $\mathsf I_2(2k+1)$ ($k\geqslant2$)}}\\\hline
\multicolumn{2}{|c|}{\mbox{type}}&r&h&\gamma&d&A,B&\alpha,\beta
&\nu\\\hline\hline
\mathsf I_2(m)&(m\geqslant4)&2&m&2m^2-5m+6&\frac m2&2m-4,\beta&m-2,\beta
&0\\\hline
\multicolumn{1}{c}{\phantom{\mathsf I_2(2k+1)}}&
\multicolumn{1}{l}{\phantom{(k\geqslant3)}}&
\multicolumn{1}{c}{\phantom{0}}&
\multicolumn{1}{c}{\phantom{2k+1}}&
\multicolumn{1}{c}{\phantom{8k^2-10k+6}}&
\multicolumn{1}{c}{\phantom{10}}&
\multicolumn{1}{c}{\phantom{r,2(r-2)}}&
\multicolumn{1}{c}{\phantom{2k-2,\beta}}&
\multicolumn{1}{c}{\phantom{r-2}}
\end{array}$$\mbox{}\\[-12mm]

The table shows that in the cases where $\beta$ has a well-defined value
(and $\alpha=m_2-1$), this value is $m_3-1$ except for $\mathsf D_r$
($r\geqslant7$), where $\beta=m_{\lfloor(r+1)/2\rfloor}-1$. With the redefinition
of $d$ and $\nu$ for the types $\mathsf I_2(2k+1)$ ($k\geqslant2$) the
formula $h=\frac d2(r+2+\nu)$ is true in general, and it is also true for
$\mathsf H_2=\mathsf I_2(5)$ with the original parameters $d=2$ and $\nu=1$.

\begin{proposition}\label{Vpmexplicit}
The equality \textup{(\ref{expsum})} in
Proposition~\textup{\ref{rootsystcyclotomic}}
holds if the multisets $V_\pm$ are given as
\begin{align*}
V_-&=\{d,\,2d-2+\nu\}\mbox{ and}\\
V_+&=\{4d-4+d\nu,\,h-d-(d-1)\nu\}
\end{align*}
with $d=\frac m2$ and $\nu=0$ for $\mathsf I_2(m)$ ($m\geqslant4$);
and for $\mathsf H_2=\mathsf I_2(5)$ the original values $d=2$ and $\nu=1$
also work.
\end{proposition}
The choice in Proposition~\ref{Vpmexplicit} of the irrelevant parameters
is thus $\alpha=\beta=1$
for type $\mathsf A_1$ and as shown in the following table.
\renewcommand{\arraystretch}{1.3}%
$$\begin{array}{|c||c|c|c|c|c|c|}\hline
\mbox{type}&\mathsf A_r&\mathsf C_r/\mathsf B_r&\mathsf G_2&
\mathsf H_2\mbox{ with $d=2$, $\nu=1$}&\mathsf H_3&\mathsf I_2(m)
\mbox{ with $d=\frac m2$, $\nu=0$}\\\hline
\beta&r&r&3&2&6&\frac m2\\\hline
\end{array}$$
\renewcommand{\arraystretch}{1.2}%
\begin{prf}
Let us first look at those exceptional types for which $d\mid h$
(including $\mathsf I_2(m)$ ($m\geqslant5$)). Here we have $\nu=0$
and the (multi)set of exponents is
$$\bigl\{m_1,\dots,m_r\bigr\}=
\Bigl\{1+jd\Bigm|0\leqslant j\leqslant\frac hd-2\Bigr\}
\cup\Bigl\{2d-1+jd\Bigm|0\leqslant j\leqslant\frac hd-2\Bigr\}$$
(see \cite[Theorem~3.2~(i)]{Bu} adding $\mathsf H_4$ and $\mathsf I_2(m)$)
so that
\begin{align*}
\sum_{i=1}^r q^{m_i}&=\sum_{j=0}^{\frac hd-2}\bigl(q^{1+jd}+q^{2d-1+jd}
\bigr)=q(1+q^{2d-2})\sum_{j=0}^{\frac hd-2}q^{jd}
=\frac{q(1-q^{4d-4})(1-q^{h-d})}{(1-q^d)(1-q^{2d-2})}
\end{align*}
in agreement with the expressions for $V_{\pm}$ (with $\nu=0$).

For the remaining types we use the following table.
$$\begin{array}{|cl|c|c|c|c|c|}\hline
\multicolumn{2}{|c|}{\mbox{type}}&h&d&\nu&4d-4+d\nu,\,h-d-(d-1)\nu&d,\,
2d-2+\nu\\\hline\hline
\mathsf A_r&(r\geqslant1)&r+1&1&r&r,r&1,r\\
\mathsf C_r/\mathsf B_r&(r\geqslant2)&2r&2&r-2&2r,r&2,r\\
\mathsf D_r&(r\geqslant4)&2r-2&2&r-4&2r-4,r&2,r-2\\
\mathsf E_7&&18&4&0&12,14&4,6\\
\mathsf H_2&&5&2&1&6,2&2,3\\
\mathsf H_3&&10&4&0&12,6&4,6\\\hline
\end{array}$$
This is in agreement with the table before Proposition~\ref{Vpmexplicit}.
\end{prf}

\begin{remark}
For the $\mathsf{DE}$ types one has $V_-=\bigl\{\frac a2,\frac b2\bigr\}$
and $V_+=\bigl\{b,\frac{ra}{4}\bigr\}$, where the parameters $a$ and $b$ are
as in Kostant's article \cite{Ko}.
Note also that for those types $\frac a2=d$ and $\frac b2=\frac{h+2}{2}-d$.
We can already look ahead and use (\ref{sumAB}) to obtain $h=dr-4d+6$;
from (\ref{prodAB}) and $h^2=\gamma$ (still for the $\mathsf{DE}$ types)
and using the equality $h=dr-4d+6$ we get $d(h-2r-6d+26)=24$.
\end{remark}

\section{Synthesis and further computations}
Proposition~\ref{rootsystcyclotomic} shows that Theorem~\ref{mainthm}
can be applied in the context of root systems with $V_+=\{A,B\}$
and $V_-=\{\alpha,\beta\}$ as in the table before Proposition~\ref{Vpmexplicit}.

Define $\gamma_0,\gamma_1,\gamma_2,\gamma_3,\ldots$ (depending on
a parameter $p$) by the series expansion
\begin{equation}\label{defgammaABab}
\sum_{n=0}^\infty\gamma_n t^n=\frac{(1-\alpha t)(1-\beta t)}{(1-A t)
(1-B t)}\sqrt[\raisebox{0.3em}{$\textstyle p$}]{\frac{1+pt}{1-pt}}\,.
\end{equation}

The series expansions
\begin{align}\label{abABfraction}
\frac{(1-\alpha t)(1-\beta t)}{(1-A t)(1-B t)}&=
\bigl(1-(\alpha+\beta)t+\alpha\beta t^2\bigr)\sum_{n=0}^\infty\Bigl(
\sum_{j=0}^nA^jB^{n-j}\Bigr)t^n
\end{align}
and
\begin{align}\label{pexpansion}
\sqrt[\raisebox{0.3em}{$\textstyle p$}]{\frac{1+pt}{1-pt}}&=
\left(\sum_{j=0}^\infty\binom{\frac1p}{j}(pt)^j\right)
\left(\sum_{k=0}^\infty\binom{-\frac1p}{\,k}(-pt)^k\right)
=:\sum_{n=0}^\infty p_nt^n\\\notag
&=1+2t+2t^2+\frac{2p^2+4}{3}t^3+\frac{4p^2+2}{3}t^4
+\frac{6p^4+20p^2+4}{15}t^5+\cdots
\intertext{specializing for $p=1$ and $p=2$}\notag
\frac{1+t}{1-t}&=1+2\sum_{n=1}^\infty t^n\\\notag
\sqrt{\frac{1+2t}{1-2t}}&=\sum_{n=0}^\infty\binom{2n}{n}(1+2t)t^{2n}
=1+2t+2t^2+4t^3+6t^4+12t^5+\cdots
\end{align}
can be used to write down an explicit formula for $\gamma_n$ defined in
(\ref{defgammaABab}).

Note that the series expansion of
$\bigl((1+pt)/(1-pt)\bigr)^{1/p}$ has
integer coefficients if $p=2^k$ with $k\in
\mathbb Z_{\geqslant0}$. In fact, for $f(t)=1+\sum_{n=1}^\infty a_nt^n$
we let
$$Tf(t):=\sqrt{f(2t)}=1+\sum_{n=1}^\infty b_nt^n.$$
A comparison of coefficients shows that
$$b_n=2^{n-1}a_n-\frac12\sum_{j=1}^{n-1}b_jb_{n-j},$$
and hence if $a_1$ is even and all $a_n$ are integers,
then all $b_n$ are even. Starting with
$f(t):=(1+t)/(1-t)=1+2\sum_{n=1}^\infty t^n$,
we get $\bigl((1+2^kt)/(1-2^kt)\bigr)^{1/2^k}=T^kf(t)\in1+2t\mathbb Z[\![t]\!]$.

(Note also that in the limit $p\to0$ we get the power series expansion
of $e^{2t}$, which is a fixed point of the transformation $T$.)

\begin{remark}
The transformation $T$ on (generating series of) integer sequences
starting with $1$ and having an even integer as next term may be investigated.
Here is a tiny list of examples:
$$\begin{array}{lcl}
a_0,a_1,a_2,\ldots&\quad\stackrel{T}{\longmapsto}\quad&b_0,b_1,b_2,\ldots\\\hline
a_n=n+1&&b_n=2^n\\
a_n=2^n&&b_n=\binom{2n}{n}\\
a_n=C_{n+1}=\frac{1}{n+2}\binom{2n+2}{n+1}&&b_n=2^n C_n\\

\end{array}$$
More generally, one may fix a positive integer $\ell$ and look at the
transformation
$$f(t)\longmapsto
\sqrt[\raisebox{0.1em}{$\scriptstyle \ell$}]{f(\ell t)}$$
for $f(t)=1+\sum_{n=1}^\infty a_nt^n$ with $\ell\mid a_1$ and $a_n\in\mathbb Z$.
\end{remark}

\begin{lemma}
The elementary symmetric polynomials in $A$ and $B$ can be written as
follows.
\begin{align}\label{sumAB}
A+B&=h-2+\alpha+\beta\\\label{prodAB}
AB&=h^2-\gamma+(h-2)(\alpha+\beta-1)
+\alpha\beta
\end{align}
Furthermore,
\begin{align}\notag
X_n&:=\sum_{j=0}^nA^jB^{n-j}\\\label{Chin}
&\phantom{:}=\displaystyle\sum_{j=0}^{\bigl\lfloor\frac n2\bigr\rfloor}(-1)^j\binom{n-j}{j}
(h-2+\alpha+\beta)^{n-2j}\bigl(h^2-\gamma+(h-2)(\alpha+\beta-1)
+\alpha\beta\bigr)^j.
\end{align}
\end{lemma}
\begin{prf}
From (\ref{defgammaABab}) we get using (\ref{gammaonetwo})
\begin{align*}
\gamma_1&=2+(A+B)-(\alpha+\beta)&&{}=h\\
\gamma_2&=2+2(A+B)+(A^2+AB+B^2)\\
&\phantom{==}{}
-2(\alpha+\beta)-(\alpha+\beta)(A+B)+\alpha\beta
&&{}=\gamma-h
\end{align*}
and solving for the elementary symmetric polynomials in
$A$ and $B$ we get (\ref{sumAB}) and (\ref{prodAB}).

Note that for $n\geqslant2$
\begin{align*}
X_n&=(A+B)\sum_{j=0}^{n-1}A^jB^{n-1-j}
-AB\sum_{j=0}^{n-2}A^jB^{n-2-j}\\
&=(h-2+\alpha+\beta)X_{n-1}-\bigl(h^2-\gamma+(h-2)(\alpha+\beta-1)
+\alpha\beta\bigr) X_{n-2}
\end{align*}
with $X_0=1$ and $X_1=h-2+\alpha+\beta$.
By solving the recursion we have (\ref{Chin}).
\end{prf}

\begin{proposition}\label{gammaintermsofhrab}
The invariant $\gamma$ is expressible as a polynomial in $h,r,\alpha,\beta$,
namely,
\begin{equation}\label{gammainhrab}
\gamma=h^2+(h-2)(\alpha+\beta-1)-(r-1)\alpha\beta.
\end{equation}
\end{proposition}
\begin{prf}
From (\ref{ABrab}) we have $AB=r\alpha\beta$. The formula
(\ref{gammainhrab}) follows by combining with (\ref{prodAB}).
\end{prf}

\begin{remark}
The formula $h=\frac d2(r+2+\nu)$ follows by inserting into $AB=r\alpha\beta$
the expressions in Proposition~\ref{Vpmexplicit} for $V_+=\{A,B\}$ and
$V_-=\{\alpha,\beta\}$ and using the fact that the product
$\bigl(d(\nu-r)+4(d-1)\bigr)(d-2)\nu$ vanishes.
\end{remark}

To summarize we state the following theorem.
\begin{theorem}
Let $m_1\leqslant\cdots\leqslant m_r$ be the exponents of an irreducible
(crystallographic (and reduced) or noncrystallographic) finite
root system (of rank $r$) with Coxeter number $h$ and parameters $\gamma$
and $d$ as in the table before Proposition~\textup{\ref{Vpmexplicit}}. Put
\begin{align}\notag
\alpha&:=\begin{cases}
\mbox{arbitrary}&\mbox{if $r=1$,}\\
m_2-1&\mbox{if $r\geqslant2$,}\\
\quad\mbox{or}\\
d
\end{cases}
\intertext{and define}\label{betafromalpha}
\beta&:=\begin{cases}
\mbox{arbitrary}&\mbox{if $h=(r-1)\alpha+2$,}\\
\dfrac{h^2-\gamma+(h-2)(\alpha-1)}{2+(r-1)\alpha-h}&
\mbox{if $h\neq(r-1)\alpha+2$.}
\end{cases}
\end{align}
Let
$$\sum_{n=0}^\infty \gamma_nt^n=
\bigl(1-(\alpha+\beta)t+\alpha\beta t^2\bigr)\Bigl(\sum_{n=0}^\infty X_nt^n
\Bigr)\Bigl(\sum_{n=0}^\infty p_nt^n\Bigr)$$
with $X_n$ as in \textup{(\ref{Chin})} and $p_n$ as in
\textup{(\ref{pexpansion})}. (So $\gamma_n$ a polynomial in
$h,\gamma,\alpha,\beta$ (or, by \textup{(\ref{gammainhrab})}, alternatively
in $h,r,\alpha,\beta$)
(symmetric in $\alpha,\beta$) and depends on an additional parameter $p$
which can be chosen arbitrarily.) Then
\begin{equation}\label{powersum}
\sum_{i=1}^r m_i^{\,n}=n!\,r\Td_n(\gamma_1,\dots,\gamma_n).
\end{equation}
\end{theorem}
\begin{prf}
As already mentioned, this is an application of Theorem~\ref{mainthm}
in the context of root systems, that is, using
Proposition~\ref{rootsystcyclotomic} with $V_+=\{A,B\}$ and
$V_-=\{\alpha,\beta\}$ and inserting (\ref{abABfraction}), (\ref{pexpansion}),
and (\ref{Chin}) into the series expansion (\ref{defgammaABab}).
The expression for $\beta$ follows from Proposition~\ref{gammaintermsofhrab}.

For $r=1$ there is nothing more to say. So let's assume $r\geqslant2$.
Here is the reason why we can take $\alpha=d$ instead of $\alpha=m_2-1$:
$d=m_2-1$ in all cases except possibly for types $\mathsf I_2(m)$,
but then we get $\beta=m_2-1=m-2$. Or still slightly more generally:
for the types $\mathsf A_r$ ($r\geqslant2$), $\mathsf C_r/\mathsf B_r$,
$\mathsf I_2(m)$, and $\mathsf H_3$ we could choose $\alpha\neq m_2-1$
and automatically get $\beta=m_2-1$ from (\ref{betafromalpha}).
\end{prf}

Let's continue by writing down $\gamma_3$ and $\gamma_4$ in terms of
$h,\gamma,\alpha,\beta$ (and $p$)
\begin{align}\label{gammathree}
\gamma_3&=-h^3+2h\gamma-2\gamma+\frac13(2p^2+4)-(h^2-\gamma-h+2)(\alpha+\beta)
-(h-2)\alpha\beta\\\notag
\gamma_4&=-h^4+h^2\gamma+\gamma^2+3h^3-6h\gamma-h^2+2\gamma
+\frac23h(p^2+5)-2
\\\notag
&\phantom{==}{}-(h^2-\gamma-h+2)\bigl((2h-2+\alpha+\beta)(\alpha+\beta)
-\alpha\beta\bigr)\phantom{\frac00}\\\label{gammafour}
&\phantom{==}{}-(h-2)(2h-2+\alpha+\beta)\alpha\beta.\phantom{\frac00}
\end{align}
By inserting (\ref{gammaonetwo}), (\ref{gammathree}), and (\ref{gammafour})
into (\ref{powersum}) using the formulae $\Td_4(c_1,c_2,c_3,c_4)=
\frac{1}{720}\bigl(-c_1^{\,4}+4c_1^{\,2}c_2+c_1c_3+3c_2^{\,2}-c_4\bigr)$ and
$\Td_5(c_1,c_2,c_3,c_4,c_5)=\frac{1}{1440}\bigl(-c_1^{\,3}c_2+3c_1c_2^{\,2}
+c_1^{\,2}c_3-c_1c_4\bigr)$ we get
\begin{align}\label{expfourthpowers}
\sum_{i=1}^rm_i^{\,4}&
=\frac{r}{30}\bigl(
-h^4+5h^2\gamma+2\gamma^2-7h^3-2h\gamma+4h^2-2\gamma-2h+2
+R_{45}\bigr)\\\label{expfifthpowers}
\sum_{i=1}^rm_i^{\,5}&
=\frac{r}{12}h\bigl(2\gamma^2-2h^3-2h\gamma+4h^2-2\gamma-2h+2
+R_{45}\bigr)
\end{align}
where
\begin{align}\label{R45}
R_{45}&=
(h^2-\gamma-h+2)\bigl((h-2+\alpha+\beta)(\alpha+\beta)-\alpha\beta\bigr)
+(h-2)(h-2+\alpha+\beta)\alpha\beta.
\end{align}

Surely, one could continue and give explicit formulae for higher power sums.
Let's stop here and display formulae for the sum of the heights cubes
and fourth powers.
\begin{proposition} With $R_{45}$ as in \textup{(\ref{R45})} above we have
\begin{align*}
\sum_{\varphi\in\Phi_+}\rootht(\varphi)^3&=
\frac{r}{120}\bigl(-h^4+5h^2\gamma+2\gamma^2-7h^3+13h\gamma-6h^2+3\gamma-7h+2
+R_{45}\bigr)\\
\sum_{\varphi\in\Phi_+}\rootht(\varphi)^4&=
\frac{r}{60}(h+1) \bigl(2\gamma^2-3h^3+3h\gamma-2\gamma-3h+2+R_{45}\bigr).
\end{align*}
\end{proposition}
\begin{prf}
Insert (\ref{lowpowers}), (\ref{expfourthpowers}), and (\ref{expfifthpowers})
into (\ref{heightspowersum}).
\end{prf}

\begin{remark} Using the power series expansions for (\ref{defgammaABab})
one computes the following explicit expressions for the quantities
$\gamma_n$.
For the types $\mathsf A_r$ one gets for $n\geqslant1$
\begin{align*}
\gamma_n|_{p=1}&=r^n+r^{n-1}
\intertext{and has
$$\sum_{i=1}^r i^n=n!\,r\Td_n\bigl(r+1,r^2+r,\dots,r^n+r^{n-1}\bigr)$$
as an alternative to Bernoulli's formula (\ref{Bernoulli}).}
\intertext{For the types $\mathsf C_r$ ($r\geqslant2$) one gets}
\gamma_n|_{p=1}&=(2r)^n-2\sum_{j=0}^{n-2}(2r)^j
\intertext{but it looks somewhat more natural to specialize to $p=2$}
\gamma_n|_{p=2}&=-2\sum_{j=0}^{\bigl\lfloor\frac{n}{2}\bigr\rfloor}C_{j-1}\,(2r)^{n-2j}
\end{align*}
where $C_k=\frac{1}{k+1}\binom{2k}{k}$ is the $k$th Catalan number (for
$k\geqslant0$) and employing
the $(-1)$st Catalan number $C_{-1}=-\frac12$.
\end{remark}

One may ask whether as an alternative to our considerations using
generating series a more geometric/combinatorial approach
via toric geometry/counting lattice points in polytopes can be found
(see also \cite[Section~2.4]{BR}, where the Bernoulli polynomials are recognized
as lattice point enumerators of certain pyramids).

\end{document}